[1994/12/01]% LaTeX date must December 1994 or later
\documentclass[10pt, reqno]{amsart}
\usepackage{a4wide}
\usepackage[english, activeacute]{babel}
\usepackage{amsmath,amsthm,amsxtra}
\usepackage{epsfig}
\usepackage{amssymb}
\usepackage{latexsym}
\usepackage{amsfonts}
\usepackage{hyperref}
\title{On the ill-posedness of the 5th-order Gardner equation}
\author{Miguel A. Alejo}
\author{Eleomar Cardoso Jr.}
\address{Universidade Federal de Santa Catarina, Brazil} %University of Copenhagen \\ 2100, Denmark}
\email{miguel.alejo@ufsc.br}
\email{eleomar.junior@ufsc.br}
\date{\today}
\subjclass[2000]{Primary 37K15, 35Q53; Secondary 35Q51, 37K10}
\keywords{Gardner equation, Breather, ill-posedness, integrability}
\chardef\bslash=`\\ % p. 424, TeXbook
\newtheorem{thm}{Theorem}[section]
\newtheorem{cor}[thm]{Corollary}

\newtheorem*{thma}{{\bf Main Theorem}}

\newtheorem{defn}[thm]{Definition}

\theoremstyle{remark}
\newtheorem{rem}{Remark}[section]

\numberwithin{equation}{section}

\newcommand{\bel}{\begin{equation}\label}
\newcommand{\eeq}{\end{equation}}
\newcommand{\R}{\mathbb{R}}

\newcommand{\al}{\alpha}
\newcommand{\bt}{\beta}

\def\bm{\left( \begin{array}{cc}}
\def\endm{\end{array}\right)}

\newcommand{\be}{\begin{equation}}
\newcommand{\ee}{\end{equation}}
\newcommand{\ba}{\begin{equation*}}
\newcommand{\ea}{\begin{equation*}}
\newcommand{\bea}{\begin{eqnarray}}
\newcommand{\eea}{\end{eqnarray}}
\newcommand{\bee}{\begin{eqnarray*}}
\newcommand{\eee}{\end{eqnarray*}}
\newcommand{\ben}{\begin{enumerate}}
\newcommand{\een}{\end{enumerate}}

\newcommand{\eval}[2][\right]{\relax
  \ifx#1\right\relax \left.\fi#2#1\rvert}

\begin{document}
\begin{abstract}
We present ill-posedness results for the initial value problem (IVP) of the 5th-order Gardner equation. We use
new breather solutions discovered for this higher order Gardner equation to measure the regularity of the
Cauchy problem in Sobolev spaces $H^s(\R)$. We find the sharp Sobolev index under which the local well-posedness of the problem is lost,
meaning that the dependence of 5th-order Gardner solutions upon initial data fails to be continuous.
%
%
% We present ill-posedness results for the initial value problem (IVP) of the 5th-order Gardner equation.
% % overcoming the problem of the loss of the scaling property of this equation.
%  We measure the regularity of the Cauchy problem in the classical Sobolev spaces $H^s$,
%  and show the critical Sobolev index under which the local well-posedness  of the problem is not present,
%  in the sense that the dependence of solutions upon initial data fails to be continuous.
\end{abstract}
\maketitle \markboth{On the ill-posedness of the 5th-order Gardner equation}{Miguel A. Alejo and Eleomar Cardoso Jr.}
\renewcommand{\sectionmark}[1]{}
%\maketitle \markboth{Nonlinear stability of higher order mKdV breathers}{Miguel A. Alejo and Eleomar Cardoso}
%\renewcommand{\sectionmark}[1]{}
%%%%%%%%%%%%%%%%%%%%%%%%%%%%%%%%%%%%%%%%%%%%%%%%%%%%%%%%%%%%%%%%%%%%%%%%%%%%%%%%%%%%%%%%%%%%%%%%%%%%%%%%%%%%%%%%%%%%%%%%%%%%%%%%%%%%%%%%%%%%%%%%%%%%%%%%%%%%%%%%%%%%%%%%%%%%%%%%%%%%%%%%%%%%%%%%%%%%%%%%%%%%%%%%%%%%%%%%%%%%%%%%%%%%%%%%%%%%%%%%%%%%%%%%%%%%%%%%%%%%%%%%%%%%%%%%%%%%%%%%%%%%%%%%%%%%%%%%%%%%%%%
\section{Introduction}
In this short note we continue our work on the Gardner equation, started in \cite{Ale3}, but this time we are going to focus on the
regularity  for the Cauchy problem of the 5th-order Gardner equation (5th-GE hereafter)

%%\be\label{edoB5}\begin{aligned}
%\be\label{GE0}\begin{aligned}\begin{cases}
%%\bea\label{GE0}\begin{cases}
%&v_t + \Big(v_{4x} + 10(\mu+v)v_x^2+ 10(\mu+v)^2v_{xx}\\
%& +  60\mu^3v^2 + 60\mu^2v^3 + 30\mu v^4 + 6v^5 \Big)_x =0,~~~\mu, %t, x\in\mathbb{R},\\
%% & v_t+ \Big(v_{4x} + 10(\mu+v)v_x^2+ 10(\mu+v)^2v_{xx}\\
%% &+ 30\mu^4v +  60\mu^3v^2 + 60\mu^2v^3 + 30\mu v^4 + 6v^5 \Big)_x %=0,~~~\mu, t, x\in\mathbb{R},\\
%%+ 6v^2(10\mu^3 + 10\mu^2v + 5\mu v^2 + v^3)  \Big)_x =0,~~~\mu, t, %x\in\mathbb{R},\\
%& v(0,x)=v_0(x),
%\end{cases}\end{aligned}\ee

%\medskip
%\noindent where $v = v(t,x)$ is a real valued function. We obtain %$$\eqref{GE0} from the 5th-order mKdV equation \eqref{5th_mKdV}

%\be\label{edoB5}\begin{aligned}
\be\label{GE0}\begin{aligned}\begin{cases}
%\bea\label{GE0}\begin{cases}
&v_t + 10\mu^2v_{xxx}+v_{5x}+\displaystyle\left[\mathcal{K}_{\mu}(v)\right]_x =0,\ ~~~\mu, t, x\in\mathbb{R},\\
% & v_t+ \Big(v_{4x} + 10(\mu+v)v_x^2+ 10(\mu+v)^2v_{xx}\\
% &+ 30\mu^4v +  60\mu^3v^2 + 60\mu^2v^3 + 30\mu v^4 + 6v^5 \Big)_x =0,~~~\mu, t, x\in\mathbb{R},\\
%+ 6v^2(10\mu^3 + 10\mu^2v + 5\mu v^2 + v^3)  \Big)_x =0,~~~\mu, t, x\in\mathbb{R},\\
& v(0,x)=v_0(x),
\end{cases}\end{aligned}\ee

%%\be\label{edoB5}\begin{aligned}
%\be\label{GE0}\begin{aligned}\begin{cases}
%%\bea\label{GE0}\begin{cases}
%&v_t + \Big(v_{4x} + 10(\mu+v)v_x^2+ 10(\mu+v)^2v_{xx}\\
%& +  60\mu^3v^2 + 60\mu^2v^3 + 30\mu v^4 + 6v^5 \Big)_x =0,~~~\mu, t, x\in\mathbb{R},\\
%% & v_t+ \Big(v_{4x} + 10(\mu+v)v_x^2+ 10(\mu+v)^2v_{xx}\\
%% &+ 30\mu^4v +  60\mu^3v^2 + 60\mu^2v^3 + 30\mu v^4 + 6v^5 \Big)_x =0,~~~\mu, t, x\in\mathbb{R},\\
%%+ 6v^2(10\mu^3 + 10\mu^2v + 5\mu v^2 + v^3)  \Big)_x =0,~~~\mu, t, x\in\mathbb{R},\\
%& v(0,x)=v_0(x),
%\end{cases}\end{aligned}\ee

\medskip
\noindent where
\be\label{GE01}\mathcal{K}_{\mu}(v):=10(\mu+v)v_x^2+ 20\mu vv_{xx}+10v^2v_{xx}+ 60\mu^3v^2 + 60\mu^2v^3 + 30\mu v^4 + 6v^5,\ee

\medskip
\noindent and $v = v(t,x)$ is a real valued function. Note that we get \eqref{GE0} by looking for solutions of the 5th-order mKdV equation 

\begin{eqnarray}\label{5th_mKdV}u_t+\displaystyle\left(u_{4x}+10uu_x^2+10u^2u_{xx}+6u^5\right)_x=0,\end{eqnarray}
\noindent
in the form $u(t,x)=\mu + v(t,x),~\mu\in\R,$ and with a suitable spatial translation. Hence  for small $\mu\ll1,$
\eqref{GE0} can be considered as a perturbed 5th-order mKdV equation.  Due to this relation with \eqref{5th_mKdV}, the physical context where  \eqref{GE0}  appears is
mainly as a perturbed model for unidirectional propagation of shallow water waves over flat surfaces, but also in  wave interaction and elastic media.
See \cite{Kaku,KichOlv,MarchSm,Olv1} for further reading.

\medskip

In this work, our aim is to extend  previous results  \cite{Ale3} on the ill-posedness for the Gardner equation

\bea\label{GE00} & v_t + (v_{xx} + 6\mu v^2 + 2v^3)_x =0,~~~\mu, t,
x\in\mathbb{R},\
%+ 6v^2(10\mu^3 + 10\mu^2v + 5\mu v^2 + v^3)  \Big)_x =0,~~~\mu, t, x\in\mathbb{R},\\
%& v(0,x)=v_0(x),
\eea
\noindent
to the above \eqref{GE0} higher order version of the Gardner equation.
Specifically, we present here one result about the ill-posedness of the 5th-GE for given data
in $H^s(\R)$. Our interest is to study the nonlinear evolution driven by the 5th-GE \eqref{GE0} with initial
data defined by a periodic in time and spatially localized function as it was firstly introduced by Kenig, Ponce and Vega \cite{KePV}, and
as we will show further.\\

\medskip

About the well-posedness of the 5th-order mKdV equation (\ref{5th_mKdV}), Linares  by  using a contraction mapping argument
in \cite{Lin}, showed  that  the Cauchy problem for the equation (\ref{5th_mKdV}) is locally well-posed at $H^2(\R)$. Moreover, Kwon
\cite{Kwon} obtained that the 5th-order mKdV equation (\ref{5th_mKdV}) is locally well-posed at $H^s(\R)$, with $s\geq \frac{3}{4}$. About
global well-posedness of the 5th-order mKdV equation, we can see that the equation (\ref{5th_mKdV}) has this property at $H^s(\R)$, if $s\geq 1$ (see \cite{Gru},
\cite{Kwon} and \cite{Lin} for more details). On the other hand,  Alejo and Kwak \cite{AK} have obtained local and global well-posedness results
for the 5th-GE \eqref{GE0} in $H^s(\R),~s\geq2$, generalizing a previous well-posedness \cite{Ale2} for the Gardner equation \eqref{GE00}.
Finally, and as far as we know, the study of the well posedness in the periodic setting for \eqref{GE0} is actually an interesting open problem, and which in fact
it would generalize recent results of Kwak for the 5th-order mKdV \eqref{5th_mKdV} in the torus \cite{Kwak}.

\medskip

Now, we remember some basic concepts:  by local well-posedness
of the IVP \eqref{GE0} we understand that there exists a unique
solution $u(t,\cdot)$ of \eqref{GE0} taking values in $H^s$ for a
time interval $[0, T)$, it defines a continuous curve in $H^s$ and
depends continuously on the initial data, that
is:\\\\
{\it for any $\epsilon>0$, there exists $\delta>0$, such that if $||u_{01}-u_{02}||_{H^s}<\delta$, then $||u_1-u_2||_{H^s}<\epsilon$, with
$\delta = \delta(\epsilon,M)$, where $||u_{0i}||_{H^s}\leq M,~~i=1,2.$}\\

In \cite{AleCar1}, we were able to build explicit breather solutions for higher order Gardner equations.
In this note, we will use the corresponding definition for the real breather solution of the 5th-GE: %\footnote{Note that if $u(x,t)=\si+v(x,t)$ is a solution of mKdV \eqref{mKdVf}, then $v(x,t)$ is solution of GE \eqref{GE0}, next of a suitable spatial translation.}(see e.g. \cite{PeGr,Ale}):\newpage

\begin{defn}[Real breather solution of the 5th-GE]\label{5thbreatherG}

Let $\al, \bt,\in \R\backslash\{0\},~\mu\in \R^{+}\backslash\{0\}$ such that $\Delta=\al^2+\bt^2-4\mu^2>0$. The real breather solution of the
5th-GE \eqref{GE0} is given explicitly by the formula

\begin{eqnarray}\label{5thBresG}B_{5\mu}\equiv B_{\alpha,
\beta,\mu}(t,x;x_1,x_2):=2\partial_x\displaystyle\left[\arctan\displaystyle\left(\frac{G_{\mu}(t,x)}{F_{\mu}(t,x)}\right)\right],\end{eqnarray}
where
\begin{eqnarray}\label{Gmu} G_{\mu}(t,x):=
\displaystyle\frac{\beta\sqrt{\alpha^2+\beta^2}}{\alpha\sqrt{\Delta}}\sin(\alpha
y_1)-\displaystyle\frac{2\mu\beta[\cosh(\beta y_2)+\sinh(\beta
y_2)]}{\Delta},\end{eqnarray}
\begin{eqnarray}\label{Fmu}F_{\mu}(t,x):= \cosh(\beta
y_2)-\displaystyle\frac{2\mu\beta[\alpha\cos(\alpha
y_1)-\beta\sin(\alpha
y_1)]}{\alpha\sqrt{\alpha^2+\beta^2}\sqrt{\Delta}},\end{eqnarray}
with
\begin{eqnarray}\label{y1y2GE}y_1:=x+\delta_5t+x_1,\quad
y_2:=x+\gamma_5t+x_2,\end{eqnarray}

\noindent
and with velocities
\begin{eqnarray}\label{delta5}\delta_5:=-\alpha^4+10\alpha^2\beta^2-5\beta^4  +10(\alpha^2-3\beta^2)\mu^2 -30\mu^4 \end{eqnarray} and
\begin{eqnarray}\label{gamma5}\gamma_5:=-\beta^4+10\alpha^2\beta^2-5\alpha^4  +10(3\alpha^2-\beta^2)\mu^2 -30\mu^4.\end{eqnarray}
\end{defn}

\begin{rem} Note that since the functional form of this breather solution is the same as the \emph{classical} breather solution of the
Gardner equation \cite[Def.1.1]{Ale3}, we get similar properties on it. For instance,
% $\delta \neq \ga$, for all values of $\al$ and $\bt$ different from zero. This means that the variables $x+\delta t$ and $x+\ga t$ are always independent. Indeed, if $\delta = \ga$, one has from (\ref{deltagamma}) $$2(\al^2 +\bt^2) =0,$$ which means $\al=\bt=0$, a contradiction. Moreover, note
for each fixed time, the  breather of the 5th-GE \eqref{5thBresG} is a function in the Schwartz class, with zero mean
$$\int_\R B_{5\mu} =0.$$\end{rem}

Furthermore, $B_{5\mu}$ satisfies the same fourth order elliptic equation which characterizes \emph{classical} Gardner breather solutions (see \cite{AleCar1} for
further reading)

\be\label{edoB5}\begin{aligned}
B_{5\mu,4x} &+ 2 (\al^2-\bt^2)(B_{5\mu,xx} + 6\mu B_{5\mu}^2 + 2B_{5\mu}^3) + (\al^2+\bt^2)^2 B_{5\mu}  + 10B_{5\mu}^2B_{5\mu,xx}+ 10B_{5\mu} B_{5\mu,x}^2 \\
& + 6B_{5\mu}^5 + 10\mu B_{5\mu,x}^2 + 20\mu B_{5\mu} B_{5\mu,xx} + 40\mu^2B_{5\mu}^3 + 30\mu B_{5\mu}^4 =0.
\end{aligned}
\ee

\begin{rem} It is truth, from Definition \ref{5thbreatherG}, that
\begin{eqnarray}\label{5thBresG.2}B_{5\mu}\equiv B_{\alpha,
\beta,\mu}(t,x;x_1,x_2):=\displaystyle\frac{2M_{\mu}(t,x)}{N_{\mu}(t,x)}\end{eqnarray}
is a breather solution to the equation (\ref{GE0}), where
\begin{eqnarray}\label{Nmu}N_{\mu}(t,x)&=&\displaystyle\left[\cosh(\beta
y_2)-\displaystyle\frac{2\mu\beta[\alpha\cos(\alpha
y_1)-\beta\sin(\alpha
y_1)]}{\alpha\sqrt{\alpha^2+\beta^2}\sqrt{\Delta}}\right]^2\nonumber\\
&+&\displaystyle\left[\displaystyle\frac{\beta\sqrt{\alpha^2+\beta^2}}{\alpha\sqrt{\Delta}}\sin(\alpha
y_1)-\displaystyle\frac{2\mu\beta[\cosh(\beta y_2)+\sinh(\beta
y_2)]}{\Delta}\right]^2\end{eqnarray} and
\begin{eqnarray}\label{Mmu}M_{\mu}(t,x)&=&\displaystyle\left[\displaystyle\frac{\beta\sqrt{\alpha^2+\beta^2}}{\sqrt{\Delta}}\cos(\alpha
y_1)-\displaystyle\frac{2\mu\beta^2[\cosh(\beta y_2)+\sinh(\beta
y_2)]}{\Delta}\right]\nonumber\\
&\times&\displaystyle\left[\cosh(\beta
y_2)-\displaystyle\frac{2\mu\beta[\alpha\cos(\alpha
y_1)-\beta\sin(\alpha
y_1)]}{\alpha\sqrt{\alpha^2+\beta^2}\sqrt{\Delta}}\right]\\
&-&\displaystyle\left[\beta\sinh(\beta
y_2)+\displaystyle\frac{2\mu\beta[\alpha\sin(\alpha
y_1)+\beta\cos(\alpha
y_1)]}{\sqrt{\alpha^2+\beta^2}\sqrt{\Delta}}\right]\nonumber\\
&\times&\displaystyle\left[\displaystyle\frac{\beta\sqrt{\alpha^2+\beta^2}}{\alpha\sqrt{\Delta}}\sin(\alpha
y_1)-\displaystyle\frac{2\mu\beta[\cosh(\beta y_2)+\sinh(\beta
y_2)]}{\Delta}\right].\nonumber\end{eqnarray} \end{rem}

%Therefore,
We will assume here that $\al, \bt>0$, which we will call them as
the \emph{frequency} and the \emph{amplitude}  of the breather
respectively. % and $\ga_5$ will be for us the \emph{velocity} of the breather solution of the 5th-GE.
The following simplification on the
breather solution will help us to analyze the regularity of the IVP for
the 5th-GE. It is easy to see in (\ref{Nmu})-(\ref{Mmu}), that
selecting $\al$ large, such that $\bt/\al\ll1$, the Gardner breather
\eqref{5thBresG.2} reduces to

\begin{eqnarray}\label{5thBresG.3} v_{\alpha,\beta,\mu}(t,x):=B_{5\mu}(t,x)\approx
2\beta\cos(\alpha(x+\delta_5t))\textrm{sech}(\beta(x+\gamma_5t))\end{eqnarray}
or simply
\begin{eqnarray}\label{5thBresG.4} v_{\alpha,\beta,\mu}(t,x)\approx
\sqrt{2}\textrm{Re}[e^{i(\alpha(x+\delta_5
t))}Q_{\beta}(x+\gamma_5t)],\end{eqnarray} where $Q$ denotes the
solution of the nonlinear ODE

\begin{eqnarray}\label{ODE}Q''-Q+Q^3=0,\end{eqnarray} with
\begin{eqnarray}\label{Q}Q(\xi)=\sqrt{2}\textrm{sech}(\xi)\end{eqnarray} and
\begin{eqnarray}\label{Qbeta}Q_{\beta}(\xi)=\beta
Q(\beta\xi).\end{eqnarray}\

%\Bel{Aproxbreather} V_{\Al,\Bt,\Si}(X,T)\Approx 2\Bt\Cos(\Al (X+(\Al^2 -3\Bt^2) T))\Sech(\Bt (X+(3\Al^2 -\Bt^2) T)),\Eeq
%i.e. in the regime $\al$ large with respect to $\bt$, the breather of GE is equal to the breather of mKdV up to order $O(\frac{\bt}{\al})$ and it is independent of the parameter $\si$ since the characteristic terms of GE (bearing  $\si$) are also proportional to $\bt/\al$.

\section{Main Theorem}

The use of the approximation (\ref{5thBresG.4}) will be a key step
in the proof of the ill-posedness of the 5th-GE. The main result
concerning the IVP \eqref{GE0} is the following:

\begin{thma}[Ill-posedness of the 5th Gardner eq.]\label{Ill5thGE}~ \noindent
If $s<3/4$, the mapping data-solution $v_0\rightarrow v(t)$, with
$v(t)$ a solution of the IVP for the 5th-order Gardner equation
\eqref{GE0} is not uniformly continuous.\end{thma}

\medskip

\begin{rem} The proof is simple, and it reduces to first compute the distance between two initial data with different
frequencies $\al_i,~~i=1,2$ but  the same amplitude $\bt$. Next, we will measure the distance between  two
solutions at time $t=T$ and selecting large enough frecuencies $\al_i,~~i=1,2$,  we avoid the interaction of
the supports of these solutions. Finally, we obtain a lower bound to the distance of solutions at time $t=T$.
Selecting  the  frequencies $\al_i,~~i=1,2$ in a suitable way, and  since $s<3/4$, we will arrive to a contradiction
with the continuous dependence of the mapping data-solution.\end{rem}

\begin{rem} If $\mu=0$, the equation (\ref{GE0}) becomes the 5th-order mKdV equation \eqref{5th_mKdV}.
%\begin{eqnarray}\label{5th_mKdV}v_t+\displaystyle\left(v_{4x}+10vv_x^2+10v^2v_{xx}+6v^5\right)_x=0.\end{eqnarray}
Therefore, as a direct consequence, this Theorem also guarantees a result on the ill-posedness to the 5th-order mKdV equation (\ref{5th_mKdV}) which agrees
at least with the upper limit of the critical Sobolev index of the ill-posedness result for \eqref{5th_mKdV} obtained by Kwon \cite{Kwon}. Namely

\begin{cor}[Ill-posedness of the 5th-order mKdV eq.]\label{Ill5thmkdv}~ \noindent
If $s<3/4$, the mapping data-solution $u_0\rightarrow u(t)$, with $u(t)$ a solution of the IVP for the 5th-order mKdV equation
\eqref{5th_mKdV} is not uniformly continuous.\end{cor}

\end{rem}

\begin{proof}{\it (of the Main Theorem.)}

We consider the IVP for the 5th-order Gardner equation with initial data given by the breather solution \eqref{5thBresG},

\bea\label{GE03}\begin{cases}
&v_t + 10\mu^2v_{xxx}+v_{5x}+\displaystyle\left[\mathcal{K}_{\mu}(v)\right]_x =0,\\
& v(0,x)=v_{\alpha,\beta,\mu}(0,x),
\end{cases}\eea
\noindent
where $\mathcal{K}_{\mu}(v)$ is defined in (\ref{GE01}). With $\mu$ fixed, we take the parameter $\alpha$ large enough, such that $\displaystyle\frac{\beta}{\alpha}\ll 1$. Then, from \eqref{5thBresG.4}, the initial data reads
\begin{eqnarray}\label{5thBresG.5} v_{\alpha,\beta,\mu}(0,x)\approx \sqrt{2}\textrm{Re}[e^{i\alpha x}Q_{\beta}(x)],\end{eqnarray}
with $Q_{\beta}$ defined in (\ref{Qbeta}). We take

\begin{eqnarray}\label{albt}\beta=\alpha^{-2s}\ \ \textrm{and}\ \ \alpha_1,\alpha_2\sim \alpha.\end{eqnarray}

Observe that $\hat{Q}_{\beta}(\cdot)$ concentrates in the ball  $V_{\beta}(0)=\{\xi\in\mathbb{R};\ |\xi|<\beta\}$. First, we
calculate the $H^s$-norm of two different initial data for the 5th-GE in the regime with $\alpha$ large enough, such that
$\displaystyle\frac{\beta}{\alpha}\ll 1$:
\begin{eqnarray}\label{estimate}\|v_{\alpha_j,\beta,\mu}(0)\|^2_{H^s}\approx\|(1+|\xi|^2)^{s/2}\hat{Q}_{\beta}(\xi-\alpha_i)\|^2_{L^2}\approx C\alpha^{2s}\beta=C,\ \
j=1,2,\end{eqnarray}
\noindent
where $C$ denotes a constant. Second, we now compute the distance between these initial data

\begin{eqnarray}\label{estimate1}&&\|v_{\alpha_1,\beta,\mu}(0)-v_{\alpha_2,\beta,\mu}(0)\|^2_{H^s}\approx
\|(1+|\xi|^2)^{s/2}(\hat{Q}_{\beta}(\xi-\alpha_1)-\hat{Q}_{\beta}(\xi-\alpha_2))\|^2_{L^2}\nonumber\\
&&\leq
C\alpha^{2s}\|\hat{Q}_{\beta}(\xi-\alpha_1)-\hat{Q}_{\beta}(\xi-\alpha_2)\|^2_{L^2}\leq
C\alpha^{2s}\displaystyle\int_{-\infty}^{+\infty}\displaystyle\left|\displaystyle\int_{\xi-\alpha_1}^{\xi-\alpha_2}\displaystyle\frac{d}{d\rho}\hat{Q}_{\beta}(\rho)d\rho\right|^2\
d\xi\nonumber\\
&&\leq
C\alpha^{2s}\displaystyle\frac{|\alpha_1-\alpha_2|}{\beta^2}\displaystyle\int_{-\infty}^{+\infty}\displaystyle\int_{\xi-\alpha_1}^{\xi-\alpha_2}|\hat{Q}_{\beta}'(\rho)|^2\
d\rho d\xi\\
&&\leq
C\alpha^{2s}\displaystyle\frac{|\alpha_1-\alpha_2|}{\beta^2}\displaystyle\left(\displaystyle\int_{\rho+\alpha_2}^{\rho+\alpha_1}d\xi\right)\displaystyle\int_{-\infty}^{+\infty}|\hat{Q}_{\beta}'(\rho)|^2d\rho\nonumber\\
&&\leq
C\alpha^{2s}\displaystyle\frac{(\alpha_1-\alpha_2)^2}{\beta^2}\beta=C\alpha^{2s}(\alpha_1-\alpha_2)^2\alpha^{2s}=C(\alpha^{2s}(\alpha_1-\alpha_2))^2.\nonumber
\end{eqnarray}

Next, we consider the corresponding solutions
$v_{\alpha_1,\beta,\mu}(t)$ and $v_{\alpha_2,\beta,\mu}(t)$ at the
time $t=T$. We can see that
\begin{eqnarray}\label{estimate2}\|v_{\alpha_1,\beta,\mu}(T)-v_{\alpha_2,\beta,\mu}(T)\|^2_{H^s}\approx \alpha^{2s}\|v_{\alpha_1,\beta,\mu}(T) - v_{\alpha_2,\beta,\mu}(T)\|^2_{L^2}.\end{eqnarray}

From (\ref{5thBresG.4}), if $\alpha$ is large enough,
\begin{eqnarray}\label{apro0} v_{\alpha_j,\beta,\mu}(T,x)\approx
\sqrt{2}\textrm{Re}[e^{i(\alpha_j(x+\delta_5 T))}\beta
Q({\beta}(x+\gamma_5T))],\ j=1,2.\end{eqnarray} Moreover, note that
\begin{eqnarray}\label{apro1}\gamma_5=-\beta^4+10\alpha^2\beta^2-5\alpha^4  +10(3\alpha^2-\beta^2)\mu^2 -30\mu^4
%\gamma_5=-30\mu^4+10(3\alpha^2-\beta^2)\mu^2-(\beta^4-10\alpha^2\beta^2+5\alpha^4)
\sim
-5\alpha^4\end{eqnarray} and
\begin{eqnarray}\label{apro2}\alpha_1^4-\alpha_2^4\sim(\alpha_1^2+\alpha_2^2)(\alpha_1+\alpha_2)(\alpha_1-\alpha_2)\sim(\alpha_1-\alpha_2)\alpha^3.\end{eqnarray}

The information above shows that $v_{\alpha_j,\beta,\mu}(T)$,
$j=1,2$, concentrates in $V_{\beta^{-1}}(5\alpha_j^4T)$, $j=1,2$.
So, we basically have disjoint supports if
\begin{eqnarray}\label{apro4}\alpha^3(\alpha_1-\alpha_2)T\gg\beta^{-1}=\alpha^{2s}.\end{eqnarray}

Under this condition, we have that
\begin{eqnarray}\label{apro5}\|v_{\alpha_1,\beta,\mu}(T)-v_{\alpha_2,\beta,\mu}(T)\|^2_{L^2}\approx\|v_{\alpha_1,\beta,\mu}(T)\|^2_{L^2}+\|v_{\alpha_2,\beta,\mu}(T)\|^2_{L^2}\approx\beta\end{eqnarray}
and
\begin{eqnarray}\label{apro6}\|v_{\alpha_1,\beta,\mu}(T)- v_{\alpha_2,\beta,\mu}(T)\|^2_{H^s}\geq C\alpha^{2s}\beta=C.\end{eqnarray}

If we select
\begin{eqnarray}\label{select}\alpha_1=\alpha+\displaystyle\frac{\delta}{2\alpha^{2s}},\ \ \alpha_2=\alpha-\displaystyle\frac{\delta}{2\alpha^{2s}},\ \
\alpha_1-\alpha_2=\displaystyle\frac{\delta}{\alpha^{2s}},\end{eqnarray}
we have that
\begin{eqnarray}\label{select2}(\alpha^{2s}(\alpha_1-\alpha_2))^2=\delta^2\end{eqnarray}
and, from (\ref{apro4}),
\begin{eqnarray}\label{select3}\alpha^3\displaystyle\frac{\delta}{\alpha^{2s}}T\gg\alpha^{2s}.\end{eqnarray}
Finally, from (\ref{select3}),
\begin{eqnarray}\label{select4}T\gg\displaystyle\frac{\alpha^{4s-3}}{\delta}.\end{eqnarray}

Since $s<\displaystyle\frac{3}{4}$, given $\delta,T>0$, we can choose $\alpha$ so large that \eqref{select4} is still valid, and
then \eqref{apro6} does not satisfy uniform continuity and we conclude.

\end{proof}

%\subsection*{Acknowledgments} The author wishes to express his sincere thanks to the anonymous referee for
%pointing out some mistakes in a preliminar version of this work. The author is also indebted for its hospitality to
%the University of Bonn where this paper was written.


\begin{thebibliography}{99}
\small{

%\bibitem{AKNS}M.J. Ablowitz, D.J. Kaup, A.C. Newel~and~ H. Segur. {\it Method for solving the sine-Gordon equation}. Phys.Rev.Letters {\bf 30}, 1262-1264 (1973).

% \bibitem{Ale}  M. A. Alejo, \textit{Focusing mKdV breather solutions with nonvanishing boundary conditions by the Inverse Scattering Method.}
% Jour. of NonLin. Math. Phys. {\bf 19}, no. 1 (2012), 1--17.

% \bibitem{Ale1} M. A. Alejo, \textit{Geometric breathers of the mKdV equation.} Acta Applicandae Mathematicae {\bf 121}, Issue 1 (2012),
% 137-155.

\bibitem{Ale3} M.A. Alejo, \textit{On the ill-posedness of the Gardner equation.} Journal of Mathematical Analysis and Applications {\bf 396}, Issue 1 (2012),
256--260.

\bibitem{Ale2} M.A. Alejo, \textit{Well-posedness and stability results for the Gardner equation.} Nonlinear Differential Equations and Applications NoDEA {\bf 19},
no. 4 (2012), 503--520.

\bibitem{AleCar1} M.A. Alejo and E. Cardoso, \emph{Dynamics of breathers in the Gardner hierarchy: universality of the variational characterization,} preprint, arXiv:1901.10409.

%M.A. Alejo and E. Cardoso, \emph{On the variational structure of breather solutions III: the Gardner hierarchy}, in preparation.

%\bibitem{AleGV}    M.A. Alejo,~C. Gorria ~and~ L. Vega, \textit{Discrete Conservation laws and the convergence of long time simulations of the mKdV equation. } submitted.

\bibitem{AK} M.A. Alejo and C. Kwak, {\it The Initial Value Problem for the 5th order Gardner equation},  submitted, arXiv:1901.03350.

% \bibitem{AMV}    M. A. Alejo, C. Mu\~noz and L. Vega, \textit{The Gardner equation and the $L^2$-stability of the $N$-soliton solution of the Korteweg-de Vries equation.}
% Transactions of the American Math. Soc. \textbf{365}, no. 1 (2013),
% 195-–212.

%\bibitem{AncoNgWill}    S.C. Anco, N.T. Ngatat ~and~ M. Willouboughby, \textit{Interaction properties of complex mKdV solitons.} Physica D 240 (2011) 1378--1394.

%\bibitem{AuYF}     T.  Au-Yeung, P. C. W. Fung ~and~ C. Au,  \textit{ Modified KdV solitons with non-zero vacuum parameter
%obtainable from the ZS-AKNS inverse method}, J. Phys. A: Math. Gen. 17 (1984) 1425--1436.

%\bibitem{B}    J. Bourgain, \textit{ Fourier Transform Restriction Phenomena for Certain Lattice subsets and Applications to Nonlinear Evolution Equations,} Geometric and Functional Analysis, \textbf{3}, (1993), 107--156,~209--262.

% \bibitem{BKPSV}    B. Birnir, C. Kenig, G. Ponce, N. Svanstedt and L. Vega, {\it On the ill-posedness of the IVP for the generalized Korteweg-de Vries and nonlinear
% Schr\"odinger equations}. J. London Math. Soc. (2) \textbf{53}
% (1996), 551--559.
%
% \bibitem{BPS}    B. Birnir, G. Ponce and N. Svanstedt, {\it The local ill-posedness of the modified KdV equation}. Ann. Inst. H. Poincar\'e Anal. Non Lineaire {\bf 13}
% (1996), 529--535.


\bibitem{Gru} A. Gr\"unrock, \emph{On the hierarchies of higher order mKdV and KdV equations}. Cent. Eur. J. Math. Vol. \textbf{8} (2010) (3), 500--536.

\bibitem{Kaku} T. Kakutani, \emph{Weakly nonlinear Hydromagnetic waves in a cold collision free plasma}, J.Phys.Soc.Japan, {\bf 26} 5 1969.


\bibitem{KePV} C. Kenig, G. Ponce and L. Vega, \textit{On the Ill-posedness of some Canonical Dispersive Equations}. Duke Mathematical Journal {\bf 106}, no. 3 (2001), 617--633.
\bibitem{KichOlv} S. Kichenassamy and P.J. Olver, {\it Existence and non-existence of solitary wave solutions to higher order model evolution equations,} SIAM J. Math. Anal. {\bf 3}
(1992), 1141-1166.


%\bibitem{KePV2} C. Kenig, G. Ponce and L. Vega, {\it Well-posedness and scattering results for the generalized Korteweg-de Vries equation via the contraction principle}. Comm. Pure Appl. Math. {\bf 46}, no.4 (1993), 527--620.

%\bibitem{KePV3} {\bf-----------}, {\it The Cauchy problem for the KdV equation in Sobolev spaces of negative indices}, Duke Math. Journ. {\bf71} n.1 (1993) 1--20.

%\bibitem{KePV4} {\bf-----------}, {\it A Bilinear Estimate with Applications to the KdV equation,} Journ.Amer.Math.Soc.  {\bf8} n.2 (1996) 573--603.

%\bibitem{La}     G.L. Lamb, \textit{Elements of Soliton Theory}, Pure Appl.Math., Wiley, New York, 1980.

\bibitem{Kwak} C. Kwak, \emph{ Low regularity Cauchy problem for the fifth-order modified KdV equations on $\mathbb{T}$}, Jour. Hyperb. Diff. Equations {\bf 15},
No. 3 (2018) 463--557.
%http://dx.doi.org/10.1142/S0219891618500170
%{\it Local well-posedness for the fifth-order KdV equations on $\mathbb{T}$}, Journal of Differential Equations {\bf260}, 10, (2016), 7683-7737.

\bibitem{Kwon} S. Kwon, \emph{Well posedness and Ill-posedness of the Fifth-order modified KdV equation}. Electr. Journal Diff. Equations. {\bf 2008}, no. 1, 1--15.

\bibitem{Lin} F. Linares, \emph{A higher order modified Korteweg-de Vries equation}. Comp. Appl. Math. \textbf{14}, no. 3 (1995),  253--267.

\bibitem{MarchSm} T.R. Marchant and N.F. Smyth, {\it The extended Korteweg-de Vries equation and the resonant flow of a fluid over topography}, J. Fluid Mech. {\bf 221} 1990, 263-288.

%\bibitem{MoSaTz} L. Molinet, J.C. Saut~and~ N. Tzvetkov, {\it Well-posedness and ill-posedness results for the Kadomtsev-Petviashvili-I equation,} Duke Math. J. {\bf 115} (2002) 353--384.

%\bibitem{Mu}     C. Mu\~noz, \textit{ The Gardner equation and the stability of the multi-kink solutions of the mKdV equation.} Submitted. Preprint on arXiv:1106.0648v3.

%\bibitem{OW}     K. Ohkuma~and~M. Wadati, \textit{Multiple-Pole Solutions of the Modified Korteweg-de Vries Equation,} J.Phys.Soc.Japan, \textbf{51}, n.6, (1982), 2029--2035.
\bibitem{Olv1} P.J. Olver, {\it Hamiltonian perturbation theory and water waves,} Contemp. Math. 28 (1984), 231-249.

%\bibitem{PeGr} D. Pelinovsky and R. Grimshaw, \textit{ Structural transformation of eigenvalues for a perturbed algebraic soliton potential}. Phys. Lett. A \textbf{229}, no. 3 (1997), 165--172.

%\bibitem{T}    T. Tao, \textit{Multilinear Weighted Convolution of $L^2$ Functions, and Applications to non-linear Dispersive Equations,} Amer.J.Math. {\bf123} (2001) 839--908.

%\bibitem{Tz}  N. Tzvetkov, \textit{Remark on the local ill-posedness for KdV equation}. C. R. Acad. Sci. Paris, t. 329, Serie I (1999), 1043--1047.

%\bibitem{W}      M. Wadati, \textit{The modified Korteweg-de Vries Equation,} J.Phys.Soc.Japan, \textbf{34}, n.5, (1973), 1289--1296.

%\bibitem{Z} P.E. Zhidkov, {\it Korteweg--de Vries and Nonlinear Schr\"odinger Equations: Qualitative Theory.} Lecture Notes in Mathematics, vol. 1756. Springer, Berlin
%(2001)
%
}
\end{thebibliography}
\end{document}